\documentclass[10pt,oneside,leqno]{amsart}
\usepackage{amsxtra}
\usepackage{amsopn}
\usepackage{amsmath,amsthm,amssymb}
\usepackage{amscd}
\usepackage{amsfonts}
\usepackage{latexsym}
\usepackage{verbatim}

\theoremstyle{plain}
\newtheorem{theorem}{Theorem}[section]

\newtheorem{lemma}[theorem]{Lemma}
\newtheorem{prop}[theorem]{Proposition}
\newtheorem{cor}[theorem]{Corollary}
\newtheorem{rem}[theorem]{Remark}

\newcommand\C{{\mathbb C}}
\newcommand\R{{\mathbb R}}
\newcommand\Z{{\mathbb Z}}

\newcommand\T{{\mathbb T}} \newcommand\Proj{{\mathbb P}}
\begin{document}
\title{Non-K\"ahler solvmanifolds with generalized K\"ahler structure}
\author{Anna Fino and Adriano Tomassini }
\date{\today}
\address{Dipartimento di Matematica \\ Universit\`a di Torino\\
Via Carlo Alberto 10\\
10123 Torino\\ Italy} \email{annamaria.fino@unito.it}
\address{Dipartimento di Matematica\\ Universit\`a di Parma\\ Via G.P. Usberti 53/A\\
43100 Parma\\ Italy} \email{adriano.tomassini@unipr.it}
\thanks{This work was supported by the Projects MIUR ``Geometric Properties of Real and Complex Manifolds'',  ``Riemannian Metrics
and Differentiable Manifolds" and by GNSAGA
of INdAM}
\begin{abstract}
We construct a compact  6-dimensional solvmanifold
endowed with a  non-trivial invariant  generalized K\"ahler structure and which does
not admit any K\"ahler metric. This is  in contrast with the case of nilmanifolds which
cannot admit any  invariant generalized K\"ahler structure unless they
are tori.
\end{abstract}
\maketitle
\section{Introduction}
The generalized K\"ahler structures were introduced and studied by M. Gualtieri in his
PhD thesis \cite{Gu} in the more general context of generalized
geometry started by N. Hitchin in \cite{Hi}.

There are many explicit constructions of non-trivial generalized-K\"ahler structures
\cite{AGG,AG,Hi2,Ko,LT,BCG,DM}.  For instance Gualtieri proved that all compact-even dimensional semisimple Lie groups are generalized
K\"ahler. In  \cite{LT} the generalized K\"ahler quotient construction is considered in relation with the hyperk\"ahler quotient construction  and generalized K\"ahler structures are given on $\C \Proj^n$, on some toric varieties  and on the complex Grassmannian.

Some obstructions and conditions on the underlying complex manifolds were found (see \cite{AG,Ca,Gu} and related references).

By \cite{Gu,AG}   it turns out that a generalized K\"ahler
structure  on a  $2n$-dimensional manifold $M$ is equivalent to a pair of Hermitian structures $(J_+, g)$ and
$(J_-, g)$, where $J_\pm$ are two integrable almost complex
structures on $M$ and $g$ is
a Hermitian  metric with respect to  $J_\pm$,  such that  the $3$-form
$ H = d^c_+  F_+  = -  d^c_- F_ -$ is closed, where $F_{\pm} ( \cdot, \cdot) = g( J_{\pm} \cdot, \cdot)$ are the
fundamental $2$-forms associated with the Hermitian structures
$(J_{\pm}, g)$ and $d^c_{\pm} = i (\overline \partial_ {\pm} -
\partial_{\pm})$. In particular, any K\"ahler metric $(J, g)$ gives
rise to a generalized K\" ahler structure by taking $J_+ = J$ and $J_-
= \pm J\,.$\newline
In the context of Hermitian geometry,  the closed $3$-form $H $ is called the
{\it torsion} of the generalized K\"ahler structure  and  it can be also identified with the torsion of the
Bismut connection associated with the Hermitian structure $(J_{\pm},
g)$ (see  \cite{Bi,Ga2}). The  generalized K\"ahler structure is called {\em untwisted} or {\em twisted} according to the fact that
the cohomology class $[H] \in H^3 (M, \R)$ vanishes or not.
In this paper we will give a homogeneous example of twisted generalized K\"ahler manifold which does not admit any K\"ahler structure.

If $(J_{\pm}, g)$ is a generalized K\"ahler structure,  then the fundamental $2$-forms $F_{\pm}$ are $\partial_{\pm} \overline
\partial_{\pm}$-closed. Therefore the  Hermitian structures $(J_{\pm}, g)$ are {\em strong K\"ahler with torsion} (SKT). Such
structures have been studied by many authors in \cite{Eg,FPS,GGP,IP,SSTV}.

In dimension four a Hermitian metric $g$ which satisfies the SKT condition is {\em standard} in the terminology of  Gauduchon (\cite{Ga})
and if, in addition,  $M$ is compact, then any Hermitian conformal
class contains a standard metric.
\newline
Compact examples in six dimensions are given  in \cite{FPS} where the Hermitian manifolds are provided
by {\em nilmanifolds}, i.e. compact quotients of
nilpotent Lie groups endowed with an Hermitian structure $(J, g)$ in
which $J$ and $g$ arise from left-invariant tensors. In \cite{Ca} it
was proved that these
manifolds  cannot admit an invariant  generalized K\"ahler structure,
since they are not formal unless they are tori. Nevertheless, all
$6$-dimensional nilmanifolds admit invariant generalized complex structures \cite{CG}.

No general restrictions  are known in the case of {\em solvmanifolds},
i.e. on compact quotients of solvable Lie groups  by uniform discrete
subgroups. By \cite{Ha2} a solvmanifold has a K\" ahler structure if
and only it is covered by a complex torus which has a structure of a
complex torus bundle over a complex torus.

As far as we know, the only known solvmanifold carrying  a
generalized K\" ahler structure is the Inoue surface. In \cite{AG}
the complex solvmanifold  from \cite{DT} was considered  and it was shown that this manifold does not admit any
left-invariant SKT metric compatible with the natural left-invariant complex structure.

Other examples are given in \cite{Eg}, where the SKT metrics are called {\em pluriclosed} in its terminology. By \cite{FG}
there are compact complex manifolds of dimension
higher than four which do no admit any SKT Hermitian metric.

In \cite{AG} generalized K\"ahler structures for which the corresponding complex structures $J_{\pm}$ commute are studied
and a classification
of compact $4$-dimensional endowed with a generalized K\"ahler structure for which the induced complex structures yield opposite
orientations is obtained.  By \cite{AG} if $J_+$ and $J_-$ commute, then the product $J_+ J_-$ is an involution of  $TM$ and the tangent bundle
splits as $TM = T_- M \oplus T_+ M$ direct sum of the
$\pm 1$-eigenspaces of the involution.
By their result,  the holomorphic tangent bundle of a compact complex surface $(M^4, J)$ admitting a generalized K\" ahler structure splits as a
direct sum of two holomorphic sub-bundles and $(M^4, J)$ is biholomorphic to one of the following:
\begin{enumerate}
\item a geometrically ruled surface;
\item a bi-elliptic complex surface;
\item a compact complex surface of Kodaira dimension 1 and even first Betti number;
\item a compact complex surface of general type;
\item a Hopf surface;
\item an Inoue surface in the family constructed in \cite{In}.
\end{enumerate}

Moreover, in relation to the distinction between untwisted and
twisted generalized K\"ahler structures, they prove that untwisted
generalized K\"ahler structures on compact $4$-dimensional manifolds may exist only if the first Betti
number is even and in the twisted case the first Betti number must be odd.

The Inoue surface $S^0$ can be viewed also as the quotient
${\mathcal H} \times \C /G_{M}$, where ${\mathcal H} $ is the
upper-half of the complex plane $\C$ and $G_M$ is a group of
analytic automorphisms of ${\mathcal H}  \times \C$ (see
\cite{In,Tr}).  The complex surface $S^0$ can be also obtained as
a $4$-dimensional solvmanifold \cite{Ha}. Moreover, since by
Hattori's theorem \cite{hat}  its de Rham cohomology is given by the invariant one, it is easy
to check that $S^0$ is formal.

Generalizing the description of the Inoue surface $S^0$ as
solvmanifold,  in section 3 we will construct a  compact
$6$-dimensional manifold with  a  twisted generalized K\"ahler
structure. The $6$-dimensional  manifold is a solvmanifold,  a
compact quotient of a non-completely solvable Lie group and the
generalized Kahler structure is left-invariant. Such a manifold does
not admit any K\"ahler structure since its first Betti number is  one
and it is a total space of a $\T^2$-bundle over the Inoue
surface. The construction can be extended in any even dimension
bigger than six.

Moreover, in the last section  we give an example of  a non-unimodular $6$-dimensional Lie group endowed
with a generalized K\"ahler structure. The corresponding Hermitian structures $(J_{\pm}, g)$ are locally conformal K\"ahler.
\medskip

\noindent{\em{Acknowledgements}} We would like to thank Vestislav
Apostolov,  Gil Cavalcanti, Sergio Console and Simon Salamon  for useful comments and conversations.
\section{Preliminaries}
In this section we briefly recall the definition of generalized  complex and K\"ahler structures, following  \cite{Gu,Hi,AG,Ca}.
Let $M$ be a $m$-dimensional manifold. Denote by $TM$ and $T^*
M$ the tangent and cotangent bundle of $M$ respectively.

Let $H$ be a closed $3$-form on $M$; the (twisted)  {\em Courant bracket}  on  the sections of $TM \oplus T^* M$ is defined by
$$
\left[ \left[   X + \xi, Y + \eta \right] \right] = [X, Y] + {\mathcal L}_X \eta - {\mathcal L}_Y \xi - \frac 12 d (\eta (X) - \xi(Y)) + \iota_Y \iota_X H,
$$
where $ {\mathcal L}_X$ and $\iota_X$ denote respectively the Lie derivative and the contraction by $X$.

On $TM \oplus T^* M$ a natural symmetric pairing of signature $(m, m)$ is given by:
\begin{equation} \label{innerproduct}
< X + \xi, Y + \eta > = \frac 12 (\eta (X) + \xi (Y)).
\end{equation}

A {\em generalized complex structure} $\mathcal J$ on $(M, H)$ is a complex structure on the bundle $TM \oplus T^* M$ which preserves the pairing and whose  $i$-eigenspace  $L$ is involutive with respect to the Courant bracket.

A generalized complex structure ${\mathcal J}$ can be also viewed as  an element of  the orthogonal Lie algebra ${\mathfrak {so}}
(TM \oplus T^*M)$ and thus, with respect to the splitting
$$
{\mathfrak {so}} (TM \oplus T^*M) = \Lambda^2 TM \oplus {\mbox {End}} (TM) \oplus \Lambda^2 T^*M,
$$ it can be written as the block matrix
$$
{\mathcal J} = \left(  \begin{array} {cc}  A & \pi\\   \sigma & A
\end{array} \right),
$$
where $\pi$ is a bivector field, $A$ an endomorphism of $TM$  and $\sigma$ a $2$-form.

For  $H =0$ examples of generalized complex structures  are given by
complex and symplectic structures. In the case of a complex manifold
the $i$-eigenspace is given by $L = T^{0,1} M \oplus T^ {* 1,0}M$ and
in the symplectic  case $L =  \{ X - i \,  \iota_X \omega , X \in T_{\C} M \}$, where $\omega$ is the symplectic form on $M$. In the block matrix form we may write:
$$
{\mathcal J} _J= \left(  \begin{array} {cc}  -J  & 0\\   0 & J^*
\end{array} \right), \quad {\mathcal J}_{\omega} = \left(  \begin{array} {cc}  0 & - \omega^{-1}\\   \omega & 0
\end{array} \right).
$$

A {\em generalized K\"ahler structure}   on  a $2n$-dimensional manifold  $M$ is a pair of commuting generalized complex
structures $({\mathcal J}_1, {\mathcal J}_2)$  on $M$
which satisfy the following conditions:
\begin{enumerate}
\item[(i)] ${\mathcal J}_1$ and  ${\mathcal J}_2$ are integrable with respect to the (twisted) Courant
bracket on $TM \oplus T^* M$ and they are compatible with the natural inner product $< \cdot
, \cdot >$ of signature $(2n, 2n)$  on $TM \oplus T^* M$ given by \eqref{innerproduct};
\item[(ii)] the quadratic form $< {\mathcal J}_1 \cdot ,
{\mathcal J}_2 \cdot >$ is definite on $TM \oplus T^* M$.
\end{enumerate}

By \cite{Gu,AG} it turns out  that a generalized K\"ahler
structure  on a manifold $M$ is equivalent to a triple $(J_+, J_-, g)$, where  $J_\pm$ are two integrable almost complex
structures on $M$ and $g$ is
a Hermitian  metric with respect to  $J_\pm$,  which satisfy the equations:
\begin{equation} \label{conditionsGK}
d^c F_+ + d^c_- F_ - =0\,, \quad d (d^c_+ F_+)=0\,, \quad  d (d^c_- F_-)=0\,,
\end{equation}
where $F_{\pm} ( \cdot, \cdot) = g( J_{\pm} \cdot, \cdot)$ are the
fundamental $2$-forms associated with the Hermitian structures
$(J_{\pm}, g)$ and $d^c_{\pm} = i (\overline \partial_ {\pm} -
\partial_{\pm}) = (-1)^n J_{\pm} d J_{\pm}$. These conditions are equivalent to
\begin{equation} \label{realconditionsGK}
 J_+ d F_+ + J_- d F_- =0\,, \quad d (J_ + d F_+) =0\,, \quad d (J_ - d F_-) =0
\end{equation}
and in Physics they appear in the target space geometry for a $(2,2)$ supersymmetric
sigma model  (see e.g. \cite{GHR}).

A trivial solution of the equations (\ref{conditionsGK}) is given by a K\"ahler structure $(g, J)$ on $M$, by taking   $J_+ = J$ and $J_- = \pm J$. So the
interesting case is when $J_- \neq  \pm J_+$, i.e. when the generalized K\"ahler structure does not arise from  a K\" ahler structure.

By (\ref{conditionsGK}) the fundamental $2$-forms $F_{\pm}$ are $\partial_{\pm} \overline
\partial_{\pm}$-closed.  In general, a Hermitian structure $(J, g, F)$ is called a strong K\" ahler structure with torsion  (SKT)
if $\partial  \overline \partial F =0$ and a K\"ahler structure satisfies this condition.
\section{Compact Example}

In this section we will describe explicitly a compact 6-dimensional example of generalized K\"ahler manifold.

Consider the $2$-step solvable Lie algebra  ${\mathfrak s}_{a,b}$ with structure equations:
\begin{equation}\label{equazionidistruttura}
\left\{
\begin{array} {l}
d e^1 = a \, e^1 \wedge e^2\,,\\[3pt]
d e^2 =0\,,\\[3pt]
d e^3 = \frac{1}{2} a \, e^2 \wedge e^3\,,\\[3pt]
d e^4 = \frac{1}{2} a \, e^2 \wedge e^4\,,\\[3pt]
d e^5 = b \, e^2 \wedge e^6\,,\\[3pt]
d e^6 = - b \, e^2 \wedge e^5\,,
\end{array}
\right.
\end{equation}
 where $a, b$ are non-zero real numbers. Let  $S_{a,b}$ be  the simply-connected solvable Lie group with Lie algebra ${\mathfrak s}_{a,b}$ and
 $(t, x_1, x_2, x_3, x_4, x_5)$ be global  coordinates  on $\R^6$. Then the Lie group $S_{a,b}$ can be described using the following product:
\begin{equation} \label{product}
\left(
\begin{array}{l}
t\\
x_1\\
x_2\\
x_3\\
x_4\\
x_5
\end{array}
\right)
\cdot
\left(
\begin{array}{l}
t'\\
x'_1\\
x'_2\\
x'_3\\
x'_4\\
x'_5
\end{array}
\right)
=
\left(
\begin{array}{l}
t+t'\\
e^{a \, t} x'_1+ x_1\\
e^{\frac{a}{2} t} x'_2 + x_2\\
e^{\frac{a}{2} t} x'_3 + x_3\\
x'_4 \cos (b\, t) - x'_5 \sin (b \, t) +   x_4\\
x'_4 \sin (b \, t) + x'_5 \cos  (b \, t)+ x_5
\end{array}
\right)
\end{equation}
Then the $1$-forms
$$
\begin{array} {l}
e^1 = e^{-a \, t} dx_1\,, \quad e^2 = dt\,, \quad e^3 = e^{ \frac {a}{2} t} dx_2\,, \quad
e^4 = e^{ \frac {a}{2} t} dx_3\,, \\[3pt]
e^5 = \cos (b \, t) \, dx_4 + \sin (b \, t) \, dx_5\,, \quad e^6 = - \sin (b \, t) \, dx_4 + \cos (b \, t) \, dx_5\,.
\end{array}
$$
are left-invariant  on $S_{a,b}$ and they satisfy the structure equations \eqref{equazionidistruttura}.
It turns out that  $S_{a,b}$ is a unimodular semidirect product
$$
\R \ltimes_{\varphi} (\R \times \R^2 \times \R^2)\,,
$$
where  $\varphi = (\varphi_1, \varphi_2)$ is the diagonal action of $\R$ on $\R \times \R^2 \times \R^2$, given by  \eqref{product}.

We start with the following
\begin{lemma} The solvable Lie group $S_{1,\frac{\pi}{2}}$ $($corresponding to $a=1, b = \frac{\pi}{2})$
admits a compact quotient $M^6 = S_{1,\frac{\pi}{2}} / \Gamma$.
\end{lemma}
\begin{proof} In order to construct a uniform discrete subgroup of $S_{1,\frac{\pi}{2}}$, we will proceed as follows.
The $4$-dimensional  solvable Lie group $\R \ltimes_ {\varphi_1}(\R \times \R^2)$ with structure equations
$$
\left\{
\begin{array}{l}
d e^1 = e^1 \wedge e^2\,,\\[3pt]
d e^2 =0\,,\\[3pt]
d e^3 = \frac{1}{2} e^2 \wedge e^3\,,\\[3pt]
d e^4 = \frac{1}{2} e^2 \wedge e^4\,,\\
\end{array}
\right.
$$
admits a compact quotient by a uniform discrete subgroup of the form $\Gamma_1 = \Z \ltimes_{\varphi_1} \Z^3$, since it can
be identified with the Inoue surface $M^4$ of type $S^0$ \cite{In}, described
as in \cite{Tr, Ha}. More precisely,  the action $\varphi_1$ can be given
by assigning   a matrix  $\varphi_1(1) =(m_{jk}) \in SL (3, \Z)$,  with two conjugate eigenvalues  $\alpha, \overline \alpha$ and a irrational
eigenvalue
$c > 1$ such that $\vert \alpha \vert^2 c = 1$ and  considering the product  on $\R \ltimes (\R \times \C)$ defined by
$$
(t, u, z) \cdot ( t', u', z') = (t + t', c^t u' + u,  \alpha^t z' + z), \quad t, t' , u, u'  \in \R, z, z' \in \C.
$$
If we denote by $(\alpha_1, \alpha_2, \alpha_3)$ an eigenvector corresponding to $\alpha$ and by $(c_1, c_2, c_3)$
a real eigenvector corresponding to $c$, then $\Gamma_1$ is generated by
$$
\begin{array}{l}
h_0: (t, u, z) \mapsto   ( t + 1,   c u, \alpha z)\,,\\
h_j: (t,u,z) \mapsto  (t, u + c_j, z + \alpha_j)\,, \quad j = 1,2,3.
\end{array}
$$
The vectors $(c_j, \alpha_j)$, $j = 1,2,3$ are linearly indipendent over $\R$ and
\begin{equation} \label{eigenvectors}
(c \, c_j, \alpha \, \alpha_j) = \sum_{k=1} ^3 m_{j_k} (c_k, \alpha_k),
\end{equation}
for any $j = 1,2,3$.

The $3$-dimensional solvable Lie group $\R \ltimes \R^2$ with structure equations
$$
\left\{
\begin{array} {l}
d e^2 =0\,,\\[3pt]
d e^5 = 2 \pi e^2 \wedge e^6\,,\\[3pt]
d e^6 = - 2 \pi e^2 \wedge e^5\,.
\end{array}
\right.
$$
is a non-completely solvable Lie group which admits a compact quotient and the uniform discrete subgroup is of the form
$\Gamma_2 =  \Z \ltimes \Z^2$ (see \cite[Theorem 1.9]{OT} and \cite{Mi}). Indeed, the Lie group $\R \ltimes \R^2$  is the group of matrices
$$
\left( \begin{array} {cccc} \cos (2 \pi t )& \sin(2 \pi t)& 0 & x\\
-\sin (2 \pi t )& \cos(2 \pi t)& 0 & y\\
0&0&1&t\\
0&0&0&1 \end{array}
\right)
$$ and a lattice $\Gamma_2$ is generated by:
$$
\left( \begin{array} {cccc} \cos (\frac{2 \pi n}{p} )& \sin(\frac{2 \pi n}{p})& 0 & x\\
-\sin (\frac{2 \pi n}{p})& \cos(\frac{2 \pi n}{p})& 0 & y\\
0&0&1&\frac {n}{p}\\
0&0&0&1 \end{array}
\right),  \left( \begin{array} {cccc} 1& 0& 0 & u_1\\
0&1& 0 & v_1\\
0&0&1&0\\
0&0&0&1 \end{array}
\right), \left( \begin{array} {cccc} 1& 0& 0 & u_2\\
0&1& 0 & v_2\\
0&0&1&0\\
0&0&0&1 \end{array}
\right),
$$
where  $n$ is an integer, $p = 2,3,4,6$ and
$$
\det \left( \begin{array} {cc} u_1&v_1\\ u_2&v_2 \end{array} \right) \neq 0.
$$

Therefore $S_{1,\frac{\pi}{2}}$ is isomorphic to  $(\R^6 = \R \ltimes (\R \times  \C \times \C), *)$, where the  product $*$ is given by
$$
(t, u, z,  w) * (t', u', z',  w') = (t + t', c^t u' + u, \alpha^t z' + z, e^{i  \frac{\pi}{2} t} w' + w),
$$
for any $t,u,t',u' \in \R, z, w, z', w' \in \C$.

Then, a uniform discrete subgroup $\Gamma \cong \Z \ltimes(\Z^3 \times \Z^2)$ of $S_{1,\frac{\pi}{2}}$ is the group generated by
the transformations
\begin{equation} \label{generatori}
\begin{array} {l}
g_0:  (t, u, z,  w) \mapsto (t + 1, c u, \alpha z,  i w)\,,\\[3pt]
g_j:  (t, u, z,  w)\mapsto  (t, u + c_j, z + \alpha_j, w)\,, \quad j = 1,2,3,\\[3pt]
g_4: (t, u, z,  w)\mapsto  (t, u,  z , w+1)\,,\\[3pt]
g_5:(t, u, z,  w)\mapsto   (t, u,  z , w+i)\,.
\end{array}
\end{equation}
Indeed, $\Gamma$ is a closed subgroup of $S_{1,\frac{\pi}{2}}$ and the action of $\Gamma$ on $S_{1,\frac{\pi}{2}}$ is properly discontinuous and
without fixed points. The compactness of the quotient can be  checked.
\end{proof}

According  to  the Mostow structure theorem (see \cite{Mos}),  the solvmanifold $M^6$  can be fibered over $S^1$ with fiber a $5$-dimensional
torus $\T^5$, since the maximal connected nilpotent subgroup is the
abelian Lie group $N$ whose Lie algebra is  spanned by $( e_1, e_3, e_4, e_5, e_6)$, where $(e_1, \ldots, e_6)$ denotes the dual frame
of $(e^1, \ldots, e^6)$.
\begin{prop}
The compact manifold $M^6 = S_{1,\frac{\pi}{2}} / \Gamma$ is the total space of a $\T^2$-bundle over the Inoue surface and
$b_1 (M^6) =1$.
\end{prop}
\begin{proof}
Since $ \Gamma \cap (\R \ltimes (\R \times \C \times \C), *)$  is a uniform discrete subgroup of $(\R \ltimes (\R \times \C), \cdot)$, the map
$$
\begin{array}{lc}
\pi: & \R \ltimes (\R \times \C \times \C)  \to  \R \ltimes (\R \times \C)\,,\\
&   (t, u, z, w) \mapsto (t ,  u,  z)
\end{array}
$$
gives a fibration
$$
\pi: M^6 \to M^4\,,
$$
with fibre $\T^2$.

By \eqref{generatori} and \eqref{eigenvectors}  the generators of $\Gamma$ satisfy the following  relations:
$$
g_j g_k = g_k g_j\,, \quad  \forall j, k = 1, \ldots, 5.
$$
Moreover,
$$
\begin{array} {l}
[g_0,g_j]= g_0 g_j g_0^{-1} g_j^{-1}: (t, u, z, w) \mapsto  (t, u - c_j + c \, c_j, z - \alpha_j + \alpha \, \alpha_j,w)\,, \\[3pt]
\qquad\qquad\qquad\qquad\qquad\qquad\qquad\qquad\qquad\qquad\qquad\qquad\qquad j = 1,2,3\,,\\[3pt]
[g_0,g_4]= g_0 g_4 g_0^{-1} g_4^{-1}: (t, u, z, w)\mapsto  (t, u, z, w - 1 + i)\,,\\[3pt]
[g_0,g_5]=g_0 g_5 g_0^{-1} g_5^{-1}: (t, u, z, w) \mapsto  (t, u, z, w - 1 - i)
\end{array}
$$
and the other commutators $[g_j, g_k]$, for any $j, k = 1, \ldots, 5$ are trivial.
Hence
$$
\begin{array}{l}
[g_0,g_j]=g_1^{m_{j1}} g_2^{m_{j2}} g_3^{m_{j3}} g_j^ {-1}\,, \quad j = 1,2,3,\\[3pt]
[g_0,g_4]= {g_4}^{-1} g_5\,,\\[3pt]
[g_0,g_5]= {g_4}^{-1} {g_5}^{-1}\,.
\end{array}
$$
Since $\Gamma$ is  $2$-step solvable, it follows that  $[\Gamma, \Gamma]$ is a torsion-free abelian subgroup of $\Gamma$ and
the rank of $[\Gamma, \Gamma]$ is $5$. By definition (see \cite{Mos2})
$$
{\mbox {rank}} \, \Gamma = {\mbox {rank}} \, \Gamma / [\Gamma, \Gamma] +  {\mbox {rank}} \, [\Gamma, \Gamma];
$$
therefore
$$
\Gamma / [\Gamma, \Gamma] \cong \Z
$$
and consequently $b_1 (M^6) = 1$.
\end{proof}
\begin{rem}
{\rm It has  to be noted that we cannot apply the Hattori's theorem \cite{hat} to compute the de Rham cohomology of
the solvmanifold $M^6$ since the group
$S_{1, \frac{\pi}{2}}$ is non completely solvable.}
\end{rem}

As a direct consequence we obtain the following

\begin{cor}
$M^6$ does not admit any K\"ahler metric.
\end{cor}

Now we can prove  our main result:
\begin{theorem}  The compact manifold $M^6 = S_{1,\frac {\pi}{2}} / \Gamma$ carries a
left-invariant $($non-trivial$)$ twisted generalized K\"ahler structure.
\end{theorem}
\begin{proof} First of all we define the two almost complex structures $J_{\pm}\,,$ by setting
$$
\begin{array} {l}
\omega_+^1 = e^1 + i e^2\,, \quad \omega_+^2 = e^3 + i e^4\,, \quad \omega_+^3  = e^5 + i e^6\,,\\[3pt]
\omega_-^1 = e^1 - i e^2\,, \quad  \omega_-^2 = e^3 + i e^4\,, \quad \omega_-^3  = e^5 + i e^6\,.
\end{array}
$$
Then  by definition $(\omega^1_{\pm}, \omega^2_{\pm}, \omega^3_{\pm})$ are the $(1,0)$-forms associated with $J_{\pm}\,.$
The almost complex structures $J_{\pm}$ are both integrable. Indeed:
$$
\begin{array}{l}
d \omega^1_+ =  \frac{i}{2} \,  \omega^1_+ \wedge   \overline  \omega^1_+\,,\\[3pt]
d \omega^2_+ =  -\frac{i} {4} (  \omega^1_+ \wedge \omega^2_+  +  \omega^2_+   \wedge  \overline \omega^1_+  )\,,\\[3pt]
d \omega^3_+ = -  \frac{\pi} {4}  ( \omega^1_+ \wedge \omega^3_+ +  \omega^3_+  \wedge  \overline \omega^1_+  )\,,
\end{array}
$$
and
$$
\begin{array} {l}
d \omega^1_- = - \frac{i}{2} \,  \omega^1_- \wedge    \overline \omega^1_-\,,\\[3pt]
d \omega^2_- =  \frac{i} {4}  ( \omega^1_- \wedge \omega^2_- +   \omega^2_- \wedge  \overline\omega^1_-  )\,\\[3pt]
d \omega^3_- = \frac{\pi} {4}  ( \omega^1_- \wedge \omega^3_-  + \omega^3_- \wedge   \overline \omega^1_- )\,.
\end{array}
$$
Moreover it is easy to see that $J_{\pm}$ commute and $J_- \neq - J_+$.
Consider the Riemannian  metric $g$  defined by
\begin{equation} \label{metric}
g = \sum_{i = 1}^6 e^i \otimes e^i.
\end{equation}
Thus $g$ is  $J_{\pm}$-Hermitian. Denote by $F_{\pm}$ the fundamental $2$-form associated  with the Hermitian structures $(J_{\pm}, g)$;
by  a direct computation we have
$$
J_+d F_+ = - e^1 \wedge e^3 \wedge e^4 = - J_- dF_-\,.
$$
Since $e^1 \wedge e^3 \wedge e^4$ is a closed  and non-exact $3$-form,
the conditions (\ref{realconditionsGK}) are satisfied and $(J_{\pm}, g)$ define a non-trivial left-invariant
twisted generalized K\"ahler structure on $M^6 = S / \Gamma$.
\end{proof}

The metric $g$ given by \eqref{metric} is not flat since the Ricci tensor is diagonal, with the only non-vanishing component
${\mbox {Ric}} (e_2, e_2) = - \frac {3}{2}$ and the Hermitian structures $(J_{\pm}, g)$ are not locally conformally K\"ahler since
$$
d F_{\pm} = e^{2}\wedge e^3\wedge e^4\,.
$$
\begin{rem}{\rm  The previous construction can be extended in order to get a
non-trivial generalized K\"ahler  on a $\T^{2n}$-bundle over the
Inoue surface, by considering the $2$-step solvable Lie algebra
$$
\left\{
\begin{array} {l}
d e^1 = a \, e^1 \wedge e^2\,,\\[3pt]
d e^2 =0\,,\\[3pt]
d e^3 = \frac{1}{2} a \, e^2 \wedge e^3\,,\\[3pt]
d e^4 = \frac{1}{2} a \, e^2 \wedge e^4\,,\\[3pt]
d e^{2k + 3}  = b \, e^2 \wedge e^{2 k + 4}\,, \\[3pt]
d e^{2k + 4} = - b \, e^2 \wedge e^{2k + 3}\,, \quad k = 1, \ldots, n,
\end{array}
\right.
$$
with $a = 1$ and $b = \frac{\pi}{2}$.

The  two  integrable complex structures $J_{\pm}$, are given   by  setting
$$
\begin{array} {l}
\omega_+^1 = e^1 + i e^2\,, \quad \omega_+^2 = e^3 + i e^4\,, \quad
\omega_+^{k+2}  = e^{2k +3} + i e^{2k+4}\,,\\[3pt]
\omega_-^1 = e^1 - i e^2\,, \quad  \omega_-^2 = e^3 + i e^4\,, \quad \omega_-^{k+2}  = e^{2k +3} + i e^{2k+4}\,,
\end{array}
$$
as the associated $(1,0)$-forms.}
\end{rem}
\section{Non-compact homogeneous  Example}
In this section we construct a non-compact homogeneous example  of $6$-dimensional Lie group endowed with a
non-trivial generalized K\"ahler structure.

Let  $\mathfrak l$ be the Lie algebra  with structure equations:
\begin{equation}\label{strutturanoncompact}
\left\{
\begin{array} {l}
d f^1 = f^1 \wedge f^2\,,\\[3pt]
d f^2 =0\,,\\[3pt]
d f^3 = \frac{1}{2} f^2 \wedge f^3\,,\\[3pt]
d e^4 = \frac{1}{2} f^2 \wedge f^4\,,\\[3pt]
d f^5 =  \frac {1}{2} f^2 \wedge f^5\,,\\[3pt]
d f^6 = \frac  {1}{2} f^2 \wedge f^6\,,
\end{array}
\right.
\end{equation}
and  $L$ be  the simply-connected Lie group with Lie algebra $\mathfrak l$.

The Lie group $L$ can be described using the following product on $\R^6$ with  global coordinates $(t, y_1, y_2, y_3, y_4, y_5)$:
$$
\left(
\begin{array}{l}
t\\
y_1\\
y_2\\
y_3\\
y_4\\
y_5
\end{array}
\right)
\cdot
\left(
\begin{array}{l}
t'\\
y'_1\\
y'_2\\
y'_3\\
y'_4\\
y'_5
\end{array}
\right)
=
\left(
\begin{array}{l}
t + t'\\
e^{t} y'_1 + y_1\\
e^{\frac{1}{2} t} y'_2+y_2\\
e^{\frac{1}{2} t} y'_3+y_3\\
e^{\frac{1}{2} t} y'_4+y_4\\
e^{ \frac{1}{2} t} y'_5+y_5
\end{array}
\right)
$$
The $1$-forms
$$
\begin{array}{l}
f^1 = e^{-t} dx_1\,, \quad f^2 = dt\,, \quad f^3 = e^{- \frac 12 t} dy_2\,, \\[3pt]
f^4 = e^{- \frac 12 t} dy_3\,, \quad  f^5 = e^{- \frac 12 t} dy_4\,, \quad  f^6 = e^{- \frac 12 t} dy_5\,.
\end{array}
$$
are left-invariant  on $L$  and $L$ is a $2$-step completely solvable Lie group. Moreover, $L$ is  not-unimodular and
consequently by \cite{Mil} it does not admit any compact quotient.
\medskip
\begin{prop} The Lie group $L$ admits a left-invariant $($non-trivial$)$  generalized  K\"ahler structure.
\end{prop}
\begin{proof}
Consider  the two almost  complex structures $J_{\pm}\,$, whose $(1,0)$-forms are given by
$$
\begin{array} {l}
\theta_+^1 = f^1 + i f^2\,, \quad \theta_+^2 = f^3 + i f^4\,, \quad \theta_+^3  = f^5 + i f^6\,,\\[3pt]
\theta_-^1  = f^1 - i f^2\,, \quad  \theta_-^2 = f^3 + i f^4\,, \quad \theta_-^3  = f^5 + i f^6\,.
\end{array}
$$
We have
$$
\begin{array} {l}
d \theta^1_+ =  \frac{i}{2} \,  \theta^1_+ \wedge   \overline  \theta^1_+\,\\[3pt]
d \theta^2_+ =  -\frac{i} {4} (  \theta^1_+ \wedge \theta^2_+  +  \theta^2_+   \wedge  \overline \theta^1_+  )\,,\\[3pt]
d \theta^3_+ = -\frac{i} {4} (  \theta^1_+ \wedge \theta^3_+  +  \theta^3_+   \wedge  \overline \theta^1_+  )\,.
\end{array}
$$
and
$$
\begin{array} {l}
d \theta^1_- = - \frac{i}{2} \,  \theta^1_- \wedge    \overline \theta^1_-\,,\\[3pt]
d \theta^2_- =  \frac{i} {4}  ( \theta^1_- \wedge \theta^2_- +   \theta^2_- \wedge  \overline\theta^1_-  )\,,\\[3pt]
d \theta^3_- = \frac{i} {4}  ( \theta^1_- \wedge \theta^3_- +   \theta^3_- \wedge  \overline\theta^1_-  )\,.
\end{array}
$$

Then,   $J_{\pm}$ are both integrable, commute and  $J_- \neq - J_+\,$.
The Riemannian  metric
$$
g = \sum_{i = 1}^6 f^i \otimes f^i.
$$
is  $J_{\pm}$-Hermitian and
$$
J_+d F_+ = - f^1 \wedge f^3 \wedge f^4 -   f^1 \wedge f^5 \wedge f^6 = - J_- dF_-\,,
$$
where  $F_{\pm}$ the fundamental $2$-form associated
with the Hermitian structures $(J_{\pm}, g)$; since $f^1 \wedge f^3 \wedge f^4 + f^1 \wedge f^5 \wedge f^6 $ is a closed $3$-form,
the conditions (\ref{realconditionsGK}) are satisfied and $(J_{\pm}, g)$ define a non-trivial left-invariant generalized K\"ahler structure on $L$.
\end{proof}

The Ricci tensor is diagonal and given by
$$
{\mbox {Ric}} (g) =  f^1 \otimes f^1 -2 f^2 \otimes f^2  - \frac{1}{2} \sum_{j = 3}^6 f^j \otimes f^j\,,
$$
hence the metric $g$ is not flat.
Furthermore, in contrast with the previous example,  the Hermitian structures $(J_{\pm}, g)$ are  locally conformally K\"ahler since
$$
d F_{\pm} = f^2 \wedge  F_{\pm}.
$$

\end{document}